\documentclass[a4paper]{article}
\usepackage[T2A]{fontenc}
\usepackage[cp1251]{inputenc} 
\usepackage{amsmath,amssymb,amsthm,amsfonts,amscd} 
\usepackage{xcolor}

\usepackage{hyperref}
\begin{document}
\sloppy
\title{\sc The Arf--Kervaire invariant of framed manifolds as an obstruction to embeddability}
\author{\sc Peter M. Akhmetiev, Matija Cencelj and Du\v{s}an D. Repov\v{s}}
\date{}
\newtheorem{theorem}{Theorem}[section]
\newtheorem*{main*}{Main Theorem}
\newtheorem*{theorem*}{Theorem}
\newtheorem{lemma}[theorem]{Lemma}
\newtheorem{proposition}[theorem]{Proposition}
\newtheorem{corollary}[theorem]{Corollary}
\newtheorem{conjecture}[theorem]{Conjecture}
\newtheorem{problem}[theorem]{Problem}
\theoremstyle{definition}
\newtheorem{definition}[theorem]{Definition}
\newtheorem{remark}[theorem]{Remark}
\newtheorem*{remark*}{Remark}
\newtheorem*{example*}{Example}
\def\Z{{\mathbb Z}}
\def\R{{\mathbb R}}
\def\RP{{\mathbb R}\!{\rm P}}
\def\N{{\mathbb N}}
\def\C{{\bf C}}
\def\A{{\bf A}}
\def\D{{\bf D}}
\def\fr{{\operatorname{fr}}}
\def\st{{\operatorname{st}}}
\def\mod{{\operatorname{mod}\,}}
\def\cyl{{\operatorname{cyl}}}
\def\dist{{\operatorname{dist}}}
\def\sf{{\operatorname{sf}}}
\def\dim{{\operatorname{dim}}}
\def\Imm{{\mbox{\it Imm}}}
\maketitle
\begin{abstract}
We prove that no  $14$-connected
(resp. $30$-connected) 
stably parallelizable 
 manifold $N^{30}$
 (resp. $N^{62}$)
 of dimension
 $30$
 (resp. $62$) 
 with the Arf-Kervaire invariant 1 can be smoothly embedded
 into
 $\R^{36}$
 (resp. $\R^{83}$).
 \[  \] 
 \end{abstract}
\section{Closed stably parallelizable manifolds with a nontrivial Arf--Kervaire invariant}
 Let us consider a closed stably  framed $n$--dimensional manifold. Such a manifold is presented by the pair
 $(N^n,\Xi)$, where $N^n$ is a closed manifold of the dimension $\dim(N)=n$, and $\Xi$ is an isomorphism
 of bundles  $\Xi: \nu(N) = \R^k \times N^n$, where $\nu(N)$ is the $k$--dimensional normal bundle  of  $N^n$, $k > n+1 .$
 A stably parallelizable manifold is a stably framed manifold with the  forgotten stable framing.

 Suppose that $n=2^{\ell}-2=4l+2$ and that $N^n$ is $2l$--connected. Then $N^n$ is diffeomorphic
 to the connected sum of manifolds of the following three types
 (see  (Kreck, 2000) for the proof and further references):
\begin{itemize}
\item
 closed manifold $\Sigma^n$, homotopically equivalent to the standard $n$--dimensional sphere;
\item
 product of two standard spheres $S^{2l+1} \times S^{2l+1}$;
\item
 standard Arf-Kervaire manifold (constructed later).
\end{itemize}
A connected sum of two third type manifolds is diffeomorphic to
the standard sum of some number of first and second type
manifolds. The dimension  $\dim(H_{2l+1}(N^n;\Z/2))$ is always
even and equals  to $2p$, where $p$ is the number of the  summands of
second and third type. Following the
statement of the Hill, Hopkins
and Ravenel Theorem,  third type manifold can be constructed only for
  $\ell = 5,6$ and eventually $7$ (see Hill {\it et al.,} 2016).
   In  (Akhmetiev, 2022; p.3), the statement above is proved,
       except for possibly finite number of exceptional cases (there is an error in the proof of Proposition 2: the mapping of regular cobordism classes on p.34 is not a homomorphism).

\begin{definition}
A $2l$--connected manifold $N^n$, $n=4l+2,$ is said to have a 
nontrivial  Arf--Kervaire invariant if  $n=2^{\ell}-2$, $\ell =5,6$ or
$7$, and $N^n$ is diffeomorphic to the connected sum of a
 third type manifold and  some number of first and second type
manifolds.
\end{definition}
\subsubsection*{Standard Arf--Kervaire manifold}

The third type manifold is based on  manifold
 $M^{4l+2}_0$, $4l+2=n$. This manifold
was constructed by "plumbing" in  (Browder, 1972; Theorem
V.2.11). 
Manifold $M^{4l+2}_0$ is well-defined for all nonnegative
$l$, but the condition $\partial M_0^{4l+2}=S^{4l+1}$ is fulfilled
only for $4l+2 = 2, 6, 14, 30, 62$ and eventually $126$. For $4l+2
\ne 2,6,14,30,62$ and eventually $\ne 126$, the boundary $\partial
M_0^{4l+2}$  is $PL$--homeomorphic but not diffeomorphic to the
standard $(4l+1)$--dimensional sphere.  In these exceptional
dimensions, manifold $N^{4l+2}$ is defined as $M^{4l+2}_0$
with the standard $(4l+2)$-dimensional disk glued along the
boundary $\partial M_0^{4l+2}$. In the case $n=2,6,14,$ a second
type manifold is obtained. In the case $n=30, 62$ and eventually
$126$, the obtained manifold is of the third type.
A simplified proof of existence of $M^{4l+2}_0$ for $n=62$ appears
in (Lin, 2001).
In (Jones and Rees, 1978) there is a remark, that the manifold $N^{n}$ of the third
type is $PL$--embeddable into $\R^{n+2}$.

Our main result is formulated in the following theorem.

\begin{theorem}\label{main}
(1) Let $N^{30}$ (resp. $N^{62}$) be an arbitrary closed
 $14$-connected
(resp. $30$-connected) stably parallelizable manifold with a 
nontrivial Arf--Kervaire invariant. Then the product $N^{30}
\times I$ (resp. $N^{62} \times I$) with the interval $I=[0,1]$ is
not smoothly embeddable into the Euclidean space  $\R^{46}$ (resp.
$\R^{94}$), provided that the corresponding embedding is equipped
by a nondegenerate normal field of $9$--frames  (resp.
$10$--frames) on the complement of the Cartesian product of the
interval $I$ and a point $N^{30} \times I \setminus \{pt\} \times
I$ (resp. on $N^{62} \times I \setminus \{pt\} \times I$). 

(2) No
stably parallelizable manifold $N^{30}$ (resp. $N^{62}$) is
smoothly embeddable into the Euclidean space $\R^{36}$ (resp.
$\R^{83}$).
\end{theorem}
\begin{remark}
Obviously, Assertion 2 of Theorem $\ref{main}$ follows from Assertion 1.  
Indeed, the composition  $N^{30} \subset \R^{36}
\subset \R^{46}$ provides the embedding $N^{30} \times I \times
D^{9} \subset \R^{46}$, where $D^9$ is the standard disk. 
 The restriction of this embedding to the
submanifold $N^{30} \times I \subset N^{30} \times I \times D^{9}$
ensures the condition of stable parallelizability in Assertion 2.
Analogously, the composition $N^{62} \subset \R^{83} \subset
\R^{94}$ provides the embedding $N^{62} \times I \times D^{10}
\subset \R^{94}$, and the restriction of this embedding to the
submanifold $N^{62} \times I \subset N^{62} \times I \times
D^{10}$ ensures the condition of stable parallelizability in
Assertion 2. Nevertheless, we give an independent proof of
Assertion 2 since this proof is simpler than the proof of
Assertion 1.
\end{remark}
\begin{remark} 
Eccles constructed in (Eccles, 1979; Corollary 1.2) a
stably parallelizable $30$-dimensional (resp. $62$-dimensional)
manifold $N^{30}$ 
(resp. $N^{62}$)
with  the Arf--Kervaire invariant 1,
which is embeddable in $\R^{46}.$
\end{remark}
\begin{remark}
Other geometric applications of  the   Arf--Kervaire invariant,
related to the problem of embeddability, can be found in
(Randall, 1999).
\end{remark}
\begin{remark}
A generalization of the Kervaire invariant 1 problem and applications was given in (Akhmetiev, 2016).
\end{remark}
\begin{remark}
A preliminary version of this paper was posted on the arXiv (Akhmetiev {\it et al.,} 2010).
\end{remark}
\section{Cobordism group of immersions and the Arf--Kervaire invariant of stably framed manifolds}
Let us denote the cobordism group of immersions of oriented
$n$-manifolds in the codimension 1 by  $\Imm^{fr}(n,1)$. The class
of the regular cobordism of immersions   $f: N^n \looparrowright
\R^{n+1}$ represents an element of this cobordism group. The set of
these elements is equipped with an equivalence relation of
cobordance.
It follows by the Pontryagin-Thom construction  (Pontryagin, 1955), in the
form proposed in  (Wells, 1966), that the group $\Imm^{fr}(n,1)$ is
isomorphic to the stable $n$-homotopy group of spheres.
First, we describe the Pontryagin-Thom construction. A stably framed
manifold is a pair  $(N^n,\Xi)$ where $N^n$ is a smooth manifold
and $\Xi$ is a trivialization of the normal bundle $\nu_N$.
Namely, $N^n$ is diffeomorphic to a submanifold $M$ in
 $\R^{n+k}$. Then the normal bundle $\nu_N$ is isomorphic to
 the trivial normal bundle  $\nu_M$ and $\Xi$ is the chosen
 trivialization. The word "stably" means that $k>> n$
 (in fact $k \ge n+2$ suffices). 
 
 It is convenient to introduce
 the direct limit $k \to +\infty$. The Pontryagin-Thom construction
  (Pontryagin, 1959; Ch. 6), gives the map $F:S^{n+k} \to S^k$
  as a composition of the standard projection  $S^{n+k} \to
ME(k)(N^n)$ and the standard map $ME(k)(N^n) \to S^k$. Here,
$ME(k)(N^n)$ (or $M(N^n)$) denotes the Thom space of the trivial $k$--dimensional
bundle. This space is homeomorphic to the $k$--fold suspension of
 $N^n_+ = N^n \cup \{x\}$, where $x$ is a point. The base point
 $pt \in S^k$ is a regular value of the map $F$ and the
 preimage of a
 small neighbourhood of that point defines the framed manifold
$(N^n,\Xi)$, corresponding to the subspace of zero section of
$M(N^n)$. 

The homotopy class $[F] \in \Pi_n$ is well-defined.
Moreover, if  $F' : S^{n+k} \to S^n$ is another map, which is homotopic
to $F$, and the base point $pt$ is also a regular value of $F'$,
then it  can be constructed a framed manifold $(N',\Xi')$
analogously, with a framed cobordism $(W,\Xi_W)$, connecting
$(N,\Xi)$ and $(N',\Xi')$. Therefore  mapping $[F] \mapsto
[(N^n,\Xi)]$ defines an isomorphism between the stable homotopical
group of spheres and the cobordism group of stably framed
manifolds.

By the Smale-Hirsch theorem (Hirsch, 1961), a stably framed manifold
$(N^n,\Xi)$ defines an immersion $f: N^n \looparrowright \R^n$.
Immersion $f$ is not defined uniquely; if $f'$ is another
immersion, corresponding to $(N^n,\Xi)$, then $f$ and $f'$ are
regularly cobordant. If $(N'^n,\Xi')$ is cobordant to $(N^n,\Xi)$,
then the corresponding immersion $f': N'^n \looparrowright \R^n$,
is an element of the same regular cobordism class as the immersion
$f$. Mapping $[(N^n,\Xi)] \mapsto [f]$, constructed by Hirsch
and mapping $[F] \mapsto [(N^n,\Xi)]$, constructed by
Pontryagin, define the isomorphism between the cobordism group of
immersions $Imm^{fr}(n,1)$, and the stable homotopy group of
spheres $\Pi_n$, constructed by Wells.

Consider the case
 $n=4l+2$.

{\bf Definition of the $\Z/2$-quadratic form of an immersion and
its Arf-invariant.}
Let $f: N^{4l+2} \looparrowright \R^{4l+3}$ be the immersion,
representing an element in the group  $Imm^{fr}(4l+2,1)$. 
Homology group $H_{2l+1}(N^{4k+2};\Z/2)$ is denoted shortly by $H$. By
the Poincar\'e duality, the bilinear nondegenerate form
 $b: H \times H \to \Z/2$ is well-defined.

Take the map $S^{4l+2+k} \to M(N^{4l+2}_+)$, defined by the
Pontryagin--Thom construction. Obviously, the map fulfills  the
conditions of  (Browder, 1969; Theorem  1.4), hence a quadratic form $q: H \to
\Z/2$, associated with the form $b$, can be defined. Define the
Arf-invariant $Arf(H,q)$ as the equivalence class of $q$ in the
Witt group of quadratic forms (Browder, 1972; Sect. 4). It turns out that
$Arf(H,q)$ is an invariant of the regular cobordism class  of 
immersion $f$. It is said to be the Arf--Kervaire invariant of
$f$. Hence, by the Wells theorem, this invariant defines a
homomorphism of groups
$ \Theta: \Imm^{fr}(4l+2,1) \longrightarrow \Z/2,$
 see (Browder, 1969; Sect. 6) and (Kervaire and Milnor, 1963; Sect. 8).

The form $q$ can be constructed in a different way. Take a cycle
$x \in H$. By the Thom theorem, there exist a (possibly
nonoriented) manifold $X^{2l+1}$ and a map $i_X: X^{2l+1} \to
N^{4l+2}$ such that $i_{X,\ast}([X])=x,$ where $[X]$ is the
fundamental class of  manifold $X$. Because of general position,  it
can be assumed without loss of generality, that  map $f \circ i_X: X^{2l+1} \to
N^{4l+2}$ is an immersion with only finitely many transversal
self-intersections. Denote the self-intersection points of the
immersion $i_X$ by $\{x_1, \dots, x_s\}$. For each point $x_i$
there exists a neighbourhood consisting of two $2l+1$--disks
intersecting in $x_i$.  

Perform a surgery to obtain a manifold
$Y^{2l+1} \subset N^{4l+2}$ such that also $i_{Y,\ast}([Y])=x$,
where $i_Y: Y^{2l+1} \subset N^{4l+2}$ is the inclusion. To this
end, remove both disks and glue their boundaries by a 1-handle. The
idea of such a surgery in the case $l=0$ is known from (Pontryagin, 1959; Ch. 15, Theorem 22).
Take the immersion $f:
N^{4l+2}\looparrowright \R^n$ and consider the map  $j_Y = f \circ
i_Y: Y^{2l+1} \subset \R^n$. By the general position argument,
$j_Y$ is also an embedding and is equipped with  a cross-section
$\xi_Y$ of the normal bundle   $\nu_{j_Y}$. This cross-section is
defined by the oriented normal bundle of the immersion $f$ along the cycle.

The linking number of the framed embedding $(i_Y, \xi_Y)$ is denoted
by $lk(i_Y, \xi_Y)$ (it is defined as the linking number between 
$i_Y(Y)$ and its copy along $\xi_Y$). Define
\begin{eqnarray}\label{q}
q(x)=lk(i_Y,\xi_Y) \ (mod\; 2).
\end{eqnarray}
\begin{lemma}\label{l1}
Function
$ q: H \to \Z/2, $
given by $q(x)=lk(i_Y,\xi_Y) (\hbox{mod}~ 2)$, is well defined. It
coincides with the Browder function, constructed in 
(Browder, 1969; Lemma 1.2).
\end{lemma}
\begin{corollary}\label{cor}
Function $q(x)=lk(i_Y,\xi_Y) (\hbox{mod}~ 2)$ is the quadratic form, associated to  bilinear
form $b: H \times H \to \Z/2$.
\end{corollary}
\subsubsection*{Proof of Corollary $\ref{cor}$}
The proof can be found in (Browder, 1969; Theorem 1.4).\qed
\subsubsection*{Proof of Lemma $\ref{l1}$}
Consider a stably framed cobordism $(W^{4l+3},\Xi_W)$, connecting
given pairs  $(N^{4l+2},\Xi)$ and $(N_1^{4l+2},\Xi_1)$ of stably
framed manifolds. First suppose that the following conditions are
satisfied:

(1)  manifold $N_1^{4l+2}$ is  $2l$--connected;

(2)  cobordism $W$  consists of $i$--handles,  $1 \le i \le
2l$.

The construction of  cobordism $W$ is based on spherical surgery as
described in (Kervaire and Milnor, 1963; Sect. 5). In (Novikov, 1964; Sect. 1), 
spherical surgery is developed for a more general situation. It
follows from condition 2 that  mapping
$H_{2l+1}(N^{4l+2};\Z/2) \to H_{2l+1}(W^{4l+3};\Z/2)$, 
induced by the inclusion $N^{4l+2} \subset W^{4l+3}$, is an
isomorphism while  mapping $H_{2l+1}(N_1^{4l+2};\Z/2) \to
H_{2l+1}(W^{4l+3};\Z/2)$, induced by the inclusion
$N_1^{4l+2} \subset W^{4l+3}$, is an epimorphism.

Let $q'$ be the function, constructed by Browder. Under 
condition 1, the function $q'$ has the following geometric meaning
(Browder, 1969; the last paragraph of the proof of Theorem 3.2 and the
corresponding reference). By the Hurewicz theorem, an element
   $x \in H_{2l+1}(N_1^{4l+2};\Z/2)$ can be realized as a map of
spheres: $\varphi: S^{2l+1} \to N_1^{4l+2}$. Furthermore, $\varphi$
can be realized in its homotopy class by an embedding
$\varphi_0: S^{2l+1} \subset N_1^{4l+2}$. 

Consider the embedding
$I_{N_1} \circ \varphi_0: S^{2l+1} \to N_1^{4l+2} \subset
\R^{4l+2+k},$ where $I_{N_1}$ is the inclusion which parametrizes
the manifold $N^{4l+2}_1$. The embedding $I_{N_1} \circ \varphi_0$
is equipped by the normal vector field of $k$--frames. Then by the
Hirsch theorem, the immersion $I_{N_1} \circ \varphi_0$ is
regularly homotopic to the immersion into standard space
$\R^{4l+2} \subset \R^{4l+2+k}$. Hence the framing vectors are
parallel complements to the subspace of coordinate axes. This
immersion is denoted by $\bar \varphi_0: S^{2l+1} \looparrowright
\R^{4l+2}$. 

The stable Hopf invariant of the immersion $\bar
\varphi_0$ is defined as the number of transversal
self-intersection points. This number is denoted  by $q'(x)$ of the embedding $\varphi_0$. It 
depends neither on the choice of the  embedding $\varphi_0$ in the
homotopy class of $\varphi,$ nor on the choice of  the map  $\varphi$
realizing the homology class  $x$.

 Apply the Hirsch theorem to the embedding $i_Y: S^{2l+1} \subset \R^{4l+3}$,
 equipped with the cross-section  $\xi_Y$, to construct the
immersion $i'_Y: S^{2l+1} \looparrowright \R^{4l+2}$. The value of
$lk(i_Y,\xi_Y)$ in the right part of formula $(\ref{q})$ coincides
with the parity  of number of self-intersection points of the
immersion $i'_Y$. This proves that under condition 1 the function
$q$, defined
 in $(\ref{q})$, coincides with the function
$q'$,
constructed by Browder.

Now we prove Lemma $\ref{l1}$ in the general case. Let us consider
the cobordism $W$ under condition 2. Take an arbitrary element
$x \in H_{2l+1}(N^{4l+2};\Z/2)$ and an element $x_1 \in
H_{2l+1}(N_1^{4l+2};\Z/2)$ so that the homological class $x + x_1$ is
trivial in $H_{2l+1}(W^{4l+3};\Z/2)$. It
follows by (Browder, 1971; Lemma~\ref{o}), that $q'(x)=q'(x_1)$. We have proved that
$q'(x_1)=q(x_1)$. Let us prove the following equality
\begin{eqnarray}\label{qx}
q(x)=q(x_1).
\end{eqnarray}

Let the homology class $x$ be equal to the image of the
fundamental class under the embedding $i_Y: Y^{2l+1} \subset
N^{4l+2}$ and let the homology class $x_1$ be equal to the image
of the fundamental class under  the embedding $i_{Y_1}: Y_1^{2l+1}
\subset N_1^{4l+2}$. Suppose that the mapping of polyhedron $i_Z:
Z^{2l+2} \to W$, which realizes the singular boundary of homology
classes $x$, $x_1$, is represented by the submanifolds  $Y^{2l+1}$
and $Y_1^{2l+1}$.

 It is well known that the polyhedron $W$ can be
chosen to be a manifold in the complement of some subpolyhedron of
codimension 2. Consider the singular points and the
self-intersection curve of the polyhedron  $i_Z(Z^{2l+2})$. The
self-intersection curve of the polyhedron  $i_Z(Z^{2l+2})$ lies outside
the considered codimension 2 subpolyhedron  of the polyhedron
$Z^{2l+2}$. The boundary of self-intersection curve is the set of
critical points of  the map  $i_Z$, and the number of these points is
even. Modify the polyhedron $Z^{2l+2}$ on its regular part and
modify the map $i_Z$ by surgery in 1-handles in such a way that
the map $i_Z$ has no critical points.

Consider the immersion $f: N^{4l+2} \looparrowright \R^{4l+3}
\times \{0\}$, the immersion $f_1: N_1^{4l+2} \looparrowright
\R^{4l+3} \times \{1\}$ and the immersion $F: W^{4l+3}
\looparrowright \R^{4l+3} \times [0,1]$, such that its restriction
on the upper and the lower components of boundary coincides with
the immersions $f$ and $f_1$, respectively. Consider the pairs of
embeddings and corresponding normal  sections  $(i_Y: Y \subset
\R^{4l+3} \times \{0\};\xi_Y)$, $(i_{Y_1}: Y_1 \subset \R^{4l+3}
\times \{1\};\xi_{Y_1})$.
Take the pair of embedding and the
corresponding normal  section $(i_Z: Z^{2l+2} \looparrowright
\R^{4l+3} \times [0,1],\xi_Z)$, such that its restriction to both
components of boundary coincides with the pairs  $(i_Y, \xi_Y)$,
$(i_{Y_1};\xi_{Y_1})$, respectively. 

Obviously, the self-linking
numbers of boundary embeddings with given normal  sections are
equal modulo 2. Move $i_Y(Y)$ along $\xi_Y$; the obtained manifold is
denoted by $(i_Y(Y))'.$ Analogously, denote by $(i_{Y_1}(Y_1))'$
the manifold obtained from $i_{Y_1}(Y_1)$ by sliding along
$\xi_{Y_1}$ and by $(i_{Z}(Z))'$ the manifold obtained from
$i_{Z}(Z)$ by sliding along $\xi_{Z}$. The self-linking number of
framed embedding $(i_{Y_1}, \xi_{Y_1})$ is defined as the parity
of number of points of self-intersection of the manifold
$(i_{Y_1}(Y_1))'$ with the manifold $i_{Y_1}(Y_1)$ by homotoping
$(i_{Y_1}(Y_1))'$ to infinity. 

The self-linking number of
$(i_{Y_1}, \xi_{Y_1})$ is defined similarly. Both self-linking
numbers are congruent modulo 2 since $(i_{Z}(Z))'$ intersects
$i_{Z}(Z)$ in an even number of points and when homotoping
$(i_{Z}(Z))'$ to infinity, the intersection  of $(i_{Z}(Z))'(t)$
and $(i_{Z}(Z))$ is a collection of curves lying completely in the
regular part of polyhedra $Z'$ and $Z$. Therefore,  the boundary
of this 1-manifold consists of an even number of points and these
points are intersection points of two families of the boundary
polyhedra. Formula  $(\ref{qx})$  and Lemma  $\ref{l1}$
are thus proved. \qed

{\bf The cobordism group of stably skew-framed immersions.} Let
$(\varphi, \Psi_L)$ be a pair consisting of a $(2l+1)$-dimensional
closed manifold $\varphi: L^{2l+1} \looparrowright \R^{4l+2}$ and
of a skew-framing  $\Psi_L$ of the normal bundle $\nu_{\varphi}$,
i.e., an isomorphism $\Psi_L: \nu_{\varphi} = (2l+1)\kappa$, where
$\kappa$ is the orientation line bundle $L^{2l+1}$. It means that
$w_1(\kappa)=w_1(L^{2l+1})$. The cobordism relation of pairs is
the standard one.

 The set of all such pairs forms an abelian group
 $\Imm^{sf}(2l+1,2l+1)$ with respect to the operation of disjoint
 union.
The Pontryagin-Thom construction in the form of Wells can be
applied to this cobordism group. It induces the following
isomorphism
$$\Imm^{sf}(2l+1,2l+1) \equiv \Pi_{4l+2}(P_{2l+1}), $$
where $P_{2l+1} = \RP^{\infty}/\RP^{2l}$ is the skew projective
space and $\Pi_{4l+2}(P_{2l+1}) = \lim_{t \to +\infty}
\pi_{4l+2+t}(\Sigma^t P_{2l+1})$ is the stable homotopy group
of the $4l+2$-dimensional space $P_{2l+1}$
(Akhmetiev and Eccles, 2007).

{\bf The connecting homomorphism $\delta$.} Define the homomorphism
$$ \delta : Imm^{fr}(4l+2) \to Imm^{sf}(2l+1), $$
which is called the connecting homomorphism. It is a modification
of the transfer homomorphism of Kahn-Priddy (Eccles, 1981).

 Let the immersion $f: N^{4l+2} \looparrowright \R^{4k+3}$ represent
 an element in the group $\Imm^{fr}(4l+2,1)$. Construct a skew-framed
 immersion  $(\varphi, \Psi_{L})$, where $\varphi: L^{2l+1}
\looparrowright \R^{4l+2}$. Consider the immersion $I \circ f$,
where $I: \R^{4k+3} \subset \R^{6l+3}$ denotes the standard
embedding. The immersion $I \circ f$ is equipped with the standard
framing. Let  $g: N^{4l+2} \looparrowright \R^{6l+3}$ be an
immersion, obtained from $f$ by a small deformation which ensures
general position. Hence the immersion $g$ self-intersects
transversally. The double point manifold of immersion $g$ is denoted 
by $L^{2l+1}$. 

Let $h: L^{2l+1} \looparrowright \R^{6l+3}$ be the
parametrizing immersion of $L^{2l+1}$. The normal bundle $\nu_h$
of immersion  $h$ is naturally isomorphic to the bundle  $l
\varepsilon \oplus l \kappa$, where $\kappa$ is the line bundle
over $L^{2l+1}$, which is associated to the canonical 2-fold
covering of a double point manifold. By the
Hirsch theorem, there exists an immersion $h_1$ regularly homotopic
to $h$having its image in the subspace $\R^{4l+2} \subset
\R^{6l+3}$. The regular homotopy between the immersions $h$ and $h_1$
can be extended to the regular homotopy of normal bundles, hence
the direct summands of normal bundle are parallel to the
complementary coordinate axes of the subspace $\R^{4l+2} \subset
\R^{6l+3}$.

The immersion $h_1: L^{2l+1} \looparrowright \R^{4l+2}$ is
equipped with a skew-framing of the normal bundle, defined by the
bundle isomorphism  $\Psi_L: k \kappa \equiv \nu_{h_1}$. Starting
from the immersion $f$, we have constructed the skew-framed
immersion $(h_1,\Psi_{L})$. Define the element $\delta([f]) \in
Imm^{sf}(2l+1,2l+1)$ to be the  regular skew-framed cobordism
class $[(h_1,\Psi_{L})]$.

{\bf The Browder-Eccles invariant.} An alternative definition of
the Arf-Kervaire invariant of framed immersions was given by
(Eccles, 1981). Such a definition uses the characteristic numbers of
 double point manifolds and is based on a theorem of (Browder, 1969).
In his theorem, the Arf-Kervaire  invariant is constructed
 by means of the Adams spectral sequence. The following simplest version of
 the  definition  was given
 by Eccles
  in (Akhmetiev and Eccles, 1999).
\begin{definition}
Define the Browder-Eccles invariant $\bar \Theta(f)$ of a framed
immersion $f$ by the formula
$$\bar
\Theta(f, \Xi_{N}) = h \circ \delta(f, \Xi_{N})\pmod{2},$$
\end{definition}
\noindent
where $h \circ \delta(f, \Xi_{N})=h(I \circ f)$ is the number of self-intersection points of
the immersion
$I \circ f$.
\section{Cobordism groups of stably skew-framed immersions}
In this section we define new variants of cobordism groups, namely
the cobordism groups of stably framed immersions (stably
skew-framed immersions, respectively), i.e., the immersions which
are not framed in their image--Euclidean space but are framed only
in ambiental Euclidean spaces of sufficiently big dimensions. Such
a framing (skew-framing, respectively) is said to be the stable
framing (stable skew-framing, respectively). Cobordism groups of
stably framed and stably skew-framed immersions generalize {\it
intermediate cobordism groups}, as introduced in (Eccles, 1979). The
Arf-Kervaire and the Browder-Eccles invariants can be generalized
to the invariants defined on the cobordism group of stably framed
immersions.

The new invariants are called the
 {\it  twisted Arf-Kervaire invariant} and the {\it  twisted Browder-Eccles invariant}, respectively.
 The definition of
 the twisted Arf-Kervaire invariant is closely connected to the definition of the
  {\it Arf-changeable invariant} of framed immersions in the sense of (Jones and Rees, 1978).
\subsubsection*{Definition of stably framed cobordism groups $Imm^{stfr}(4l+2,2l+1)$.}
Let $(f, \Xi_N)$ be a pair, where $f: N^{4l+2} \looparrowright
\R^{6l+3}$ is an immersion in the codimension $2l+1$, $\Xi_N$ be a
stable framing of the manifold $N^{4l+2}$, i.e., a framing of the
normal bundle of the composition $I \circ f: N^{4l+2}
\looparrowright \R^{6l+3} \subset \R^{r}$, $r \ge 8l+6$.

The set of pairs described above is equipped by an equivalence
relation, which is given by the standard relation of cobordism. Up
to the cobordism relation the set of pairs generates an Abelian
group whose operation is determined by the disjoint union. This
group is denoted by $Imm^{stfr}(4l+2,2l+1)$.

\subsubsection*{Definition of stably skew-framed cobordism groups $Imm^{stsf}(2l+1,2l+1)$.}
Let $(\varphi, \Psi)$ be a pair, where $\varphi: L^{2l+1}
\looparrowright \R^{4l+2}$ is an immersion in  codimension
$2l+1$, $\Psi_L$ a stable skew-framing of the manifold
$L^{2l+1}$ in  codimension $2l+1$, i.e., a skew-framing of the
normal bundle of the composition $I \circ f: L^{2l+1}
\looparrowright \R^{4l+2} \subset \R^{r}$, $r \ge 4l+4$ with the
bundle $(2l+1)\kappa \oplus (r-2l-1) \varepsilon$, where
$\varepsilon$ is a trivial line bundle on $L^{2l+1}$, $\kappa$ is
a given line bundle over $L^{2l+1}$, which coincides with the
oriented line bundle over $L^{2l+1}$, since
$w_1(L^{2l+1})=w_1(\nu_{\varphi})=w_1((2l+1)\kappa) =
w_1(\kappa)$.

The set of pairs described above is equipped by an equivalence
relation, which is defined by the standard regular cobordism. The
set of equivalence classes generates an Abelian group whose
operation is determined by the disjoint union. This group is denoted
by $Imm^{stsf}(2l+1,2l+1)$.

\noindent
{\bf Homomorphisms $A: \Imm^{fr}(4l+2,1) \longrightarrow
\Imm^{stfr}(4l+2,2l+1)$, $B: \Imm^{stfr}(4l+2,2l+1) \to
\Imm^{fr}(4l+2,1)$.}
An arbitrary immersion
 $f: N^{4l+2} \looparrowright
\R^{4l+3}$ determines the immersion  $(I \circ f, \Xi_{N})$ in 
codimension $2l+1$, $I: \R^{4l+3} \subset \R^{6l+3}$, which is
stably framed in the ambient space $\R^{6l+3} \subset \R^n$. The
homomorphism $A$ is defined.

An arbitrary stably framed immersion $(f: N^{4l+2} \looparrowright
\R^{6l+3}, \Psi_{N})$ obviously induces the immersion into the
space $\R^r$, $r>8l+6$.

The Hirsch theorem, applied to this immersion of  codimension
$r-6l-3$, ensures the existence of a framed immersion $(\varphi ,
\Xi_{N})$, where $\varphi: N^{4l+2} \looparrowright \R^{4l+3}$.
Now define $B([(f,\Psi_{N})])=[(\varphi,\Xi_{N})] \in
Imm^{fr}(4l+2,1)$. Obviously, by construction  $B \circ A =
Id: Imm^{fr}(4l+2,1) \to Imm^{fr}(4l+2,1)$.

{\bf Twisted Arf-Kervaire invariant $\Theta^{st}:
\Imm^{stfr}(4l+2,2l+1) \longrightarrow \Z/2$.}
We generalize the Arf-Kervaire homomorphism $\Theta:
\Imm^{fr}(4l+1,1) \to \Z/2$ and define a homomorphism $\Theta^{st}:
\Imm^{stfr}(4l+2,2l+1) \to \Z/2$ provided that $2l+1 \ne 1,3,7$,
called {\it the twisted Arf-Kervaire invariant} such that the
following diagram commutes:
$$
\begin{array}{ccc}
\Imm^{fr}(4l+2,1) & \stackrel  {A}{\longrightarrow} & \Imm^{stfr}(4l+2,2l+1)  \\
\searrow  \Theta&  & \swarrow \Theta^{st} \\
 & \Z/2 &   \\
\end{array} \eqno(*)
$$
\qed

\noindent
{\bf Auxiliary homomorphism $\pi: H_{2l+1}(N^{4l+2};\Z/2) \to
\Z/2$, $2l+1 \ne 1,3,7$.}
It is known that for $2l+1 \ne 1,3,7$ there exist exactly two
stably trivial $S^{2l+1}$-dimensional vector $SO$--bundles - this
fact was applied in the proof of (Kervaire and Milnor, 1963; Lemma 8.3). One of these
bundles is trivial; we denote it by $E(2l+1)$. The other
bundle is nontrivial, and coincides with the tangent bundle
$T(S^{2l+1})$ over the sphere $S^{2l+1}$. We need a generalization
of this fact to the case of a $(2l+1)$--dimensional stably trivial
bundle over an arbitrary closed $(4l+2)$--dimensional manifold,
possibly nonoriented.

Let $M^{2l+1}$ be a closed, possibly nonoriented manifold, $\xi$
 a $SO(2l+1)$--bundle over $M^{2l+1}$ such that it is trivial as
a stable $SO$--bundle. Let $M_1^{2l+1}$ and $\xi$  be another
manifold and  a $SO(2l+1)$--bundle as above, respectively. Let
$W^{2l+2}$ be a $(2l+2)$--dimensional polyhedron, such that it is a
manifold in the exterior of the codimension 2 skeleton. Let the
polyhedron $W$ realize a homology between the fundamental classes
of the manifolds $M^{2l+1}$ and $M_1^{2l+1}$, i.e. $\partial W^{2l+2}
= M^{2l+1} \cup M_1^{2l+1}$. In addition, suppose that there
exists a stably trivial $SO(2l+1)$--bundle $\Xi$ over $W^{2l+2}$
such that the restrictions of $\Xi$ on $M^{2l+1}$ and on
$M_1^{2l+1}$ coincide with the bundles $\xi$ and $\xi_1$,
respectively.
\begin{lemma}\label{o}
For an arbitrary above described pair $(M^{2l+1},\xi)$ there exists
an obstruction $c(M,\xi) \in \Z/2$ to the trivialization of the
bundle $\xi$. Moreover, $c(M^{2l+1},\xi)=c(M^{2l+1}_1,\xi_1)$.
\end{lemma}
\begin{corollary}\label{deg}
  Let $f: M^{2l+1} \to S^{2l+1}$ be a map of a closed manifold to
  the standard sphere such that $\deg(f)=1 \pmod{2}$. Let $\xi$ be
  a stably trivial $SO(2l+1)$--bundle over $S^{2l+1}$ such that $c(S^{2l+1},\xi)=1$.
  Then $f^{\ast}(\xi)$ is a stably trivial $SO(2l+1)$--bundle over
  $M^{2l+1}$ and $c(M^{2l+1}, f^{\ast}(\xi))=1\pmod {2}.$
\end{corollary}
\subsubsection*{Proof of Corollary $\ref{deg}$.}
The corollary follows from the fact that there exists a homology
between $f_{\ast}([M])$
and $[S]$, 
where $[M]$ and $[S]$ are the $(2l+1)$-dimensional
fundamental classes of $M^{2l+1}$ and $S^{2l+1}$.\qed
\subsubsection*{Proof of Lemma $\ref{o}$.}
Let $(M^{2l+1},\xi)$ be the pair described in the preamble of the
lemma. Denote by  $M^{[2l]} \subset M^{2l+1}$ the complement of
the highest $(2l+1)$--dimensional cell in the skeleton of a cellular
decomposition of $M^{2l+1}$. By the dimensionality argument, the
restriction of $\xi$ on $M^{[2l]}$ is a trivial bundle. Hence
there exists a map $f: (M^{2l+1},M^{[2l]}) \to (S^{2l+1},pt)$,
such that $f^{\ast}(\psi)=\xi$, where $(S^{2l+1},\psi)$ is a
bundle over $S^{2l+1}$, satisfying the conditions of Lemma~\ref{o}.

 In the
case when $M^{2l+1}=S^{2l+1}$, Lemma~\ref{o} is true since in the proof of
 (Browder, 1969; Theorem 3.2), the obstruction $(S^{2l+1},\xi)$ is
constructed by the corresponding functional cohomological
operation. Define  $(M^{2l+1},\xi)=(S^{2l+1},\psi)$. Let
$(M_1^{2l+1},\xi_1)$ be the second pair described above and
$(W^{2l+2},\Xi)$  the homology, connecting $(M^{2l+1},\xi)$ and
$(M_1^{2l+1},\xi_1)$. If $c$   equals to zero for both pairs,
Lemma~\ref{o} is proved.  

Suppose that for at least one pair -- say
$(M^{2l+1},\xi)$ -- the value of the obstruction $(M^{2l+1},\xi)$ is
1. Consider the handle (cell) decomposition of the cobordism
$(W^{2l+2},M^{2l+1})$. The index of handles can be restricted to
$\le 2l$ (the dimension of handles to $2l+1$) as in the case when
$W^{2l+2}$ is a smooth manifold. Indeed, the handle (cell)
decomposition of $(W^{2l+2},M^{2l+1})$ can be chosen so that all
$(2l+2)$--dimensional cells retract to the $(2l+1)$--skeleton of the
polyhedron $W^{2l+2}$ without the $(2l+1)$--dimensional cells of the
upper base $M_1^{2l+1}$.

Define the map $F: (W^{2l+2},M^{2l+1}) \to (S^{2l+1}\times I,
S^{2l+1}\times \{0\})$ such that
\begin{eqnarray}\label{F}
F^{\ast} \pi^{\ast} (\psi) = \Xi,
\end{eqnarray}
where $\pi: S^{2l+1} \times I \to S^{2l+1}$ is a projection onto
the lower base. The map $F$ can be extended to the
$(2l)$--dimensional handles uniquely up to homotopy. The map $F$ can
be extended also to the $(2l+1)$--dimensional handles, but possibly
nonuniquely. By  assumption, $c(S^{2l+1},\psi)=1$. Hence the
extension of the map $F$ to the $(2l+1)$--dimensional handles can be
realized in a way that the condition $(\ref{F})$ is satisfied. Now,
on the upper base of the cobordism we have $f_1^{\ast}(\psi)=\xi_1$,
therefore by  definition, $c(M^{2l+1}_1,\xi_1)=1$. Lemma
$\ref{o}$ is thus proved. \qed

Let $(f: N^{4l+2} \looparrowright \R^{6l+3}, \Xi_N)$ be a pair
defining an element of the group $Imm^{stfr}(4l+2,2l+1)$.
Consider an arbitrary cycle $x \in H=H_{2l+1}(N^{4l+2};\Z/2)$. It
is represented by an embedding $i_Y: Y^{2l+1} \to N^{4l+2}$.
Denote shortly by $\xi$ the bundle $i_Y^{\ast}(\nu_f)$, where
$\nu_f$ is the normal bundle of an immersion $f$. Since $\nu_f$ is
a stably trivial bundle (because the manifold $N^{4l+2}$ is stably
framed by $\Xi$), for the pair $i_Y^{\ast}(\nu_f)$  the
obstruction $(Y^{2l+1},\xi)$ is defined, provided that $2l+1 \ne 1,3,7$.

Define the mapping
\begin{eqnarray}\label{pi}
\pi: H_{2l+1}(N^{4l+2};\Z/2) \to \Z/2
\end{eqnarray}
given by the formula $\pi(x) = c(Y^{2l+1},\xi)$, $y \in H =
H_{2l+1}(N^{4l+2};\Z/2)$. Lemma $\ref{o}$ implies that the value
of $\pi(x)$ does not depend on the choice of the manifold $Y^{2l+1}$
and on the choice of the embedding $l_x$, which realizes the given cycle
$x$. It can be easily verified that the mapping in $(\ref{pi})$ is a
homomorphism.
\begin{definition}
Let  $q: H \to \Z/2$, $H= H_{2l+1}(N^{4l+2};\Z/2),$ be the
quadratic form defined in $(\ref{q})$ for a stably framed manifold
$(N^{4l+2},\Xi)$. For $2l+1 \ne 1,3,7$ define the twisted
quadratic form $q^{tw}$ by the formula $q^{tw}=q+\pi: H \to \Z/2$.
The Arf invariant of this twisted quadratic form  defines a
homomorphism $\Theta^{st}: \Imm^{stfr}(4l+2,2l+1) \longrightarrow
\Z/2$, which is said to be the twisted Arf-Kervaire invariant.
\end{definition}

\noindent
{\bf Twisted Browder-Eccles invariant}. Define the invariant $\bar
\Theta^{st}: \Imm^{stfr}(4l+2,2l+1) \to \Z/2$, which is said to be
the twisted Browder-Eccles invariant, starting by the construction
of the homomorphism
\begin{eqnarray}\label{delta}
 \delta^{st}: Imm^{stfr}(4l+2,2l+1) \to
Imm^{stfs}(2l+1,2l+1).
\end{eqnarray}
Suppose that an element of the group $\Imm^{stfr}(4l+2,2l+1)$ is
represented by the pair $(f: N^{4l+2} \looparrowright \R^{6l+3},
\Xi)$. The double point manifold of immersion $f$ is denoted by
$L^{2l+1}$. This manifold is equipped by the parametrizing
immersion $\varphi': L^{2l+1} \looparrowright \R^{6l+3}$. The
corresponding  normal bundle $\nu_{\varphi'}$ admits (for a
sufficiently big natural $k$) a stable isomorphism
 $\Psi': \nu_{\varphi'} \oplus
k\varepsilon \oplus k\kappa = (2l+1+k)\varepsilon \oplus
(2l+1+k)\kappa$.

The stable isomorphism $\Psi'$ defines a stable isomorphism
$\Psi_L: \nu_L \oplus k\varepsilon = (2l+1)\kappa \oplus
k\varepsilon$, where by $\nu_L$ is denoted the $(2l+1)$--dimensional
normal bundle over $L$. By the Smale-Hirsch construction an
immersion $\varphi: L^{2l+1} \looparrowright \R^{4l+2}$ and a
stably skew-framing $\Psi_L$ of the normal bundle of this
immersion are defined. Define $\delta^{st}([(f,\Xi_{f})])$ to be
the element of the group $Imm^{stsf}(2l+1,2l+1)$ corresponding to the
pair $(\varphi,\Psi_L)$.

Consider the homomorphism $Imm^{stsf}(2l+1,2l+1)
\stackrel{h}{\longrightarrow} Imm^{\D_4}(0,4l+2) = \Z/2$, defined
as the parity of the  number of double points of stably skew-framed
immersions (this invariant is called the stably Hopf invariant).
The Browder-Eccles invariant $\bar \Theta^{st}$ is defined as the
composition
 $Imm^{stfr}(4l+2,2l+1)
 \stackrel{h \circ \delta^{st}}{\longrightarrow}
 Imm^{ \D_4}(0,4l+2) = \Z/2$.

\noindent
{\bf Subgroup $Imm^{stfr}(4l+2,2l+1)^{\ast} \subset
Imm^{stfr}(4l+2,2l+1)$.}
Define an auxiliary subgroup $Imm^{stfr}(4l+2,2l+1)^{\ast} \subset
Imm^{stfr}(4l+2,2l+1)$ as the complete preimage
of the group
$\delta^{-1}(Imm^{sf}(2l+1,2l+1) \subset Imm^{stsf}(2l+1,2l+1))$
by the homomorphism $(\ref{delta}),$
with the
following
 additional assumption: there exists a representative $[(f,\Psi_f)] \in Imm^{stfr}(4k+2,2k+1)^{\ast}$,
$f: M^{4k+2} \looparrowright \R^{6n+3}$, such that $M^{4k+2}$ is 
the  boundary of a stably parallelizable manifold $W^{4l+3}$,
$\partial(W^{4l+3})=M^{4k+2}$.
 
  It is convinient to introduce $Imm^{stfr}(4l+2,2l+1)^{\ast},$
  which is an equivalent geometrical definition of an element from the preimage of $\delta.$ This
definition is valid for all  natural $l\geq 4,$ but not  for
$l=0,1,3$.

The pair
 $(f: N^{4l+2}\looparrowright \R^{6l+3}, \Xi_{N})$
represents an element $x \in Imm^{stfr}(4l+2,2l+1)^{\ast}$ if the
following is true.
Take the pair $(\varphi: L^{2l+1} \looparrowright \R^{4l+2},
\Psi_L)$ representing the element $\delta^{st}(x)$; here
$L^{2l+1}$ is the double point manifold of the immersion $f$.
Consider the canonical covering $\bar L^{2l+1} \to L^{2l+1}$ of
the double point manifold $L^{2l+1},$ more details can be found in
(Adams, 1962). Let $\bar g: \bar L^{2l+1} \looparrowright N^{4l+2}$ be the
parametrizing immersion, $[\bar L] \in H=H_{2l+1}(N^{4l+2};\Z/2)$
be the cycle obtained as the image of the fundamental class of $\bar
L^{2l+1}$ by the immersion $\bar g$. Consider the value $\pi([\bar
L])$, where $\pi: H \to \Z/2$ was defined in $(\ref{pi})$.

Without losing generality we may assume that the stably framed immersion
$(f,\Xi)$ is chosen from the regular cobordism class so that the
manifold $\bar L^{2l+1}$ is connected. This goal can be achieved
through the 1-handles surgery on the double point manifold of
immersion $f$. The techniques of such a surgery was invented in (Haefliger, 1961). Then the condition $x \in Imm^{stsf}(4l+2,2l+1)^{\ast}$ is
equivalent to $\pi([\bar L^{2l+1}])=0$.

 The last condition is
equivalent to the fact that the pull-back $\bar g^{\ast}\nu_{f}$
over $\bar L^{2l+1}$ of the normal bundle $\nu_{f}$  is not only
stably trivial but it is also trivial (in the case $l=0,1,3$ the
normal bundle $\nu_f$ is always trivial). It means that the
immersion $\varphi: L^{2l+1} \looparrowright \R^{4l+2}$ is not
only stably skew-framed but it is also skew-framed. This is the
equivalent geometrical definition of the subgroup
$Imm^{stfr}(4l+2,2l+1)^{\ast}$.

For convenience, all homomorphisms which have been constructed,
include into the common commutative diagram:
{\small
$$
\begin{array}{ccccccc}
Imm^{fr}(4l+2,1) & \stackrel{B}{\longrightarrow} &
Imm^{stfr}(4l+2,2l+1)^{\ast} & \subset  & Imm^{stfr}(4l+2,2l+1) \\
&&&&\\
\downarrow  \delta &  & \downarrow \delta^{st}& & \downarrow \delta^{st}\\
&&&&\\
Imm^{sf}(2l+1,2l+1) & = & Imm^{sf}(2l+1,2l+1)  &
\subset & Imm^{stsf}(2l+1,2l+1) \\
&&&&\\
\downarrow h  & &  & & \downarrow h  \\
&&&&\\
Z/2 = Imm^{\D_4}(0,4l+2)&& = && Imm^{\D_4}(0,4l+2)\\
\end{array}
$$
}
\section{Arf-Kervaire and Browder-Eccles (twisted) homomorphisms coincide}
The following lemma is necessary for our proof of Theorem~\ref{main}.
\begin{lemma}\label{4.1}
The twisted Arf-Kervaire homomorphism
 $\Theta^{st}: \Imm^{stfr}(4l+2,2l+1) \to  \Z/2$
coincides with the twisted Browder-Eccles homomorphism
  $\bar \Theta^{st} : \Imm^{stfr}(4l+2,2l+1) \to \Z/2 $, on the
  subgroup
 $Imm^{stfr}(4l+2,2l+1)^{\ast} \subset Imm^{stfr}(4l+2,2l+1)$.
\end{lemma}
We shall derive Lemma~\ref{4.1} from the following lemma. 
\begin{lemma}\label{4.2}
(a) The homomorphism $\Theta^{st}$ and $\bar{\Theta}^{st}$ coincide on the subgroup $Im(B) \subset Imm^{stfr}(4l+2,2l+1)^{\ast} \subset Imm^{stfr}(4l+2,2l+1)$.

 (b) For an arbitrary element $x \in Imm^{fr}(4k+2,1),$ there exists a representative $f: M^{4l+2} \looparrowright \R^{4l+3}$, $[f]=x$, such that $M^{4l+2}$ is the boundary of a stably parallelizable manifold $W^{4l+3}, \partial(W^{4l+3})=M^{4l+2}$.
 \end{lemma}
{\bf Proof of Lemma~\ref{4.2}.} (a) This was proved in (Eccles, 1981), reformulated in
the required form in
(Akhmetiev and Eccles, 1999).

(b) This is a well-known observation by Nigel Ray, concerning Khan-Priddy transfer (Eccles, 1979; reference [7]). \qed

{\bf Proof of Lemma~\ref{4.1}.} Let $(f: N^{4l+2} \looparrowright
\R^{6l+3}, \Psi_N)$ be a stably framed immersion.  
 Let $\eta: W^{4l+3} \to
\R^{6l+3} \times \R^1_+$ be a generic mapping of a stably
parallelized manifold  $(W^{4l+3},\Psi_W)$ defining the boundary
of stably framed manifold $(N^{4l+2},\Psi)$ but in general, not the immersion
$f$ 
itself. 
The dimension of the critical point
manifold $\Sigma$ of $\eta$ equals to $2l+1$; this is less than
half of dimension of the manifold--preimage $W^{4l+3}$. The
critical points of  $\eta$ are of  type $\Sigma^{1,0}$.

By the Moren theorem (Arnold {\it et al.,} 1982; Ch. 1, Par. 9,
Sect. 6, case $k=2$), the critical point manifold has the
normal form called the extended  Whitney umbrella. The formulae
describing the singularities of Whitney umbrella $\R^{s} \to
\R^{2s-1}$ can be found in (Pontryagin, 1955; Ch. 1, Par. 4). The notion
"extended" means the inclusion of the standard singularity of
umbrella into the identical polyparametrical collection of maps.

One may assume without loss of generality, after a corresponding repair of the
singularity of map $\eta$, that the critical point manifold
$\Sigma$ satisfies the following properties.

(1) $\Sigma$ is connected with connected canonical double covering
$\bar \Sigma$.

(2) $\eta(\Sigma)$ belongs to the hyperspace $\R^{6l+3} \times
\{1\}$ and the double point manifold $K^{2l+2}$ of map $\eta$ with
the boundary $\partial K^{2l+2} = L^{2l+1} \cup \Sigma^{2l+1}$ is
regular in a small neighborhood of the boundary  with respect to
the height function on $\R^1_+$ in such way that the subspace
$\R^{6l+3} \times \{1\}$ is higher than the manifold (i.e.
$K^{2l+2}$ immerses into the subspace
$[+\varepsilon,1-\varepsilon] \times \R^{6l+3}$ outside its
regular neighborhood.

Let $\eta_{1-\varepsilon}: N_{1-\varepsilon}^{4l+2}
\looparrowright \R^{6l+3} \times \{1 - \varepsilon \}$ be the
immersion defined as the restriction of $\eta$ on the complete
preimage of the hyperspace $\R^{6l+3} \times \{1 - \varepsilon
\}$. Let $L^{2l+1}_{1 - \varepsilon}$ be the component of the
double point manifold $\eta_{1-\varepsilon}
(N_{1-\varepsilon}^{4l+2})$ in the neighborhood of the critical
point boundary $\Sigma^{2l+1}$ of $K^{2l+2}$. From the assumption
$\pi(\xi)=0$ we may deduce that the normal bundle
$\nu_{L_{1-\varepsilon}}$ of the manifold
$L^{2l+1}_{1-\varepsilon}$ is decomposed into the direct sum of
the trivial bundle $\nu_{\varepsilon} = (2l+1) \varepsilon$ and a
nontrivial bundle $\nu_{\kappa} = \nu_{\varepsilon} \otimes
\kappa$, where $\kappa$ is the orientation line bundle over
$L^{2l+1}_{1-\varepsilon}$. Since the canonical covering is
connected, $\kappa$ is nontrivial.
\subsubsection*{Construction of the stably framed immersion.}
Let us construct the stably framed immersion
 \begin{eqnarray}\label{44}
  (\xi_0: N_0^{4l+2} \looparrowright \R^{6l+3}, \Psi_0)
 \end{eqnarray}
such that the double point manifold $L^{2l+1}_0$ (equipped with a
skew-framing $\Xi_0$ of the normal bundle) coincides with an
arbitrary given skew-framed immersion.

Let us start the construction by the description of standard
immersion $g_0: S^{2l+1} \looparrowright \R^{4l+2}$ with the
self-intersection points at the origin of the coordinate system
$\R^{2l+1}_1 \oplus \R^{2l+1}_2 = \R^{4l+2}$. Let
$\R^{2l+1}_{diag}$, $\R^{2l+1}_{antidiag}$ be two coordinate
subspaces defined by means of the sum and the difference of the
base vectors in the standard coordinate spaces $\R^{2l+1}_1,
\R^{2l+1}_2$.

Consider two standard unit disks $D^{2l+1}_{1} \subset
\R^{2l+1}_1$, $D^{2l+1}_2 \subset \R^{2l+1}_2$. Take a manifold
$C$ diffeomorphic to the cylinder $S^{2l} \times I$ defined as the
collection of all the segments such that each connects a pair of
points in $\partial D^{2l+1}_1$ and $\partial D^{2l+1}_2$ with
equal coordinates. The union of two disks $D^{2l+1}_1 \cup
D^{2l+1}_2$ with $C$ (after the identification of corresponding
components of the boundary) is the image of the standard sphere
$S^{2l+1}$ by a PL-immersion $g_0$ into $\R^{4l+2}$, with one
self-intersection point at the origin. After the smoothing of
corners along $\partial C$ we obtain the  smooth immersion of
sphere under construction. 

Let us describe the manifold $N^{4l+2}_0$, the stable framing
$\Xi_{N_0}$ over this manifold and the immersion $f_0: N^{4l+2}_0
\looparrowright \R^{6l+3}$, $[(f_0,\Xi_{N_0})] \in
Imm^{stsf}(4l+2,2l+1)$. Take the embedding $\eta_0: L^{2l+1}_0
\subset \R^{6l+3}$ with the normal bundle $\nu_{L_0} = \nu_1
\oplus \nu_1 \otimes \kappa$, where $\nu_1$ is a trivial
$(2l+1)$--dimensional bundle (with the prescribed trivialization)
and $\kappa$ is the orientation line bundle over $L^{2l+1}_0$.
Consider the $(2l+2)$--dimensional bundle $\nu_1 \otimes \kappa
\oplus \varepsilon$ over $L_0^{2l+1}$ and define the manifold
$N_0^{4l+2}$ as the boundary  $S(\nu_1 \otimes \kappa \oplus
\varepsilon)$ of the disk bundle of this vector bundle.

The locally trivial fibration $p: N^{4l+2}_0 \to L^{2l+1}_0$ is
well-defined. Because $\nu_1 \otimes \kappa \oplus \varepsilon$ is
the normal bundle of $L^{2l+1}_0$, the manifold $N^{4l+2}_0$
admits an embedding in codimension 1. This embedding determines
the framing $\Psi_{N_0}$ over $N_0^{4l+2}$ such that the
constructed framed manifold $(N_0^{4l+2},\Psi_{N_0})$ is bounding.
 
Let us define the immersion $f_0: N^{4l+2}_0 \looparrowright
\R^{6l+3}$. Take the normal bundle $\nu_1 \otimes \kappa \oplus
\nu_1$ of the immersion $\eta_0$ and consider the collection of
the standard immersed spheres $g_0(S^{2l+1})$ constructed above in
each fiber of $\eta_0$. The pair $(f_0,\Psi_0)$ is the stably
framed immersion under construction.
\subsubsection*{Calculation of invariants of the constructed stably framed immersion.}
This section is devoted to the calculation of the twisted
Arf-Kervaire and the twisted Browder-Eccles invariants for the
stably framed immersion ($\ref{44}$). The group
$H=H_{2l+1}(N_0^{4l+2};\Z/2)$ is generated by two elements. The
first generator $x \in H$ is represented by a spherical fiber of
the fibration $p: N^{4l+2}_0 \to L^{2l+1}_0$. The fibration $p$
has a standard section $p^{-1}$ constructed by the trivial direct
summand in the bundle $\nu_1 \otimes \kappa \oplus \varepsilon$.
The image of the fundamental class of the base $L^{2l+1}_0$
induced by $p^{-1}$, represents the second generator $y \in H$.
Let us prove under the assumption $2l+1 \ne 1,3,7$ that the
homomorphism $\pi: H \to \Z/2$ is defined by $\pi(x)=1, \pi(y)=0$.

The condition $\pi(x)=1$ holds since for an arbitrary immersion of
a sphere $f_0: S^{2l+1} \looparrowright \R^{4l+2}$ with one
self-intersection point the corresponding  normal
$2l+1$--dimensional bundle is nontrivial.
Let us prove
that the following holds
\begin{eqnarray}\label{4}
\pi(y)=0.
\end{eqnarray}
The cycle $y \in H$ is represented by the image of the fundamental
class induced by the map of section $p^{-1}(L_0^{2l+1}) \to
N_0^{4l+2}$ of the fibration $p$. The collection of the bases in
the fibers of the subbundle $\nu_1 \subset \nu_{L_0}$ defines the
trivialization of the normal bundle of the immersion $f_0$ over
the submanifold $p^{-1}(L_0^{2l+1}) \subset N_0^{4l+2}$. This
proves $(\ref{4})$.

The twisted Arf-Kervaire invariant of a stably framed immersion
$(\xi_0,\Psi_0)$ is equal to $q(y),$ i.e., coincides with the
twisted Browder-Eccles invariant. This gives the required
computations.

Let us finish the proof of Lemma~\ref{4.1}. Consider a stably framed
immersion $(f_{1-\varepsilon}, \Psi_{N_{1-\varepsilon}})$ with the
skew-framed double point manifold $(L^{2l+1}_{1-\varepsilon},
\Xi_{L_{1-\varepsilon}})$. 
The stably
framed immersion $(\eta_{1-\varepsilon},
\Psi_{N_{1-\varepsilon}})$ is regularly cobordant to the immersion
$(f,\Psi_N)$ from the beginning of the proof of Lemma~\ref{4.1}.
Obviously, by considering the normal form of Whitney umbrella,
 $[(f_{1-\varepsilon}, \Psi_{1-\varepsilon})] \in
Imm^{stfr}(4l+2,2l+1)^{\ast} \subset Imm^{stfr}(4l+2,2l+1)$ and
the stably framed immersion $(f_{1-\varepsilon},
\Psi_{N_{1-\varepsilon}})$ is in fact, framed and the stably
skew-framed immersion $(L^{2l+1}_{1-\varepsilon},
\Xi_{L_{1-\varepsilon}})$ is in fact, skew-framed.

The normal bundle $\bar{\nu}$ of the submanifold $\bar{L}^{2l+1}_{1-\varepsilon} \subset N^{4l+2}_{1-\varepsilon}$ is
a stably trivial bundle.
Let us prove that this bundle is trivial. Indeed, the fibers of $\bar{\nu}$ over antipodal pair of points of the projection $p: \bar{L}^{2l+1}_{1-\varepsilon} \to L^{2l+1}_{1-\varepsilon}$ are cannonicaly isomorphic (by the normal form arguments) and the bundle $\bar{\nu}$ is the pull-back of a stable-trivial bundle $\nu$ over $L^{2l+1}_{1-\varepsilon}$.
This is true since the obstruction $\pi$, given by  formula~(\ref{pi}), is trivial.
By this fact, the normal bundle over $L^{2l+1}_{1-\varepsilon}$ splits into a framed subbundle $\nu_{0}$ (which in general, does not coincide with $\nu$) and a skew-framed subbundle $\nu_0 \otimes \kappa_p$, like at the beginning of the construction of the stably framed immersion
$ (\xi_0: N_0^{4l+2} \looparrowright \R^{6l+3}, \Psi_0)$
 from formula $(\ref{44}).$

Let us apply the construction $(\ref{44})$, where the standard
stably framed immersion $(f_0,\Psi_{N_0})$ is obtained so that
its self-intersection points are opposite to the points of the
skew-framed immersion
$(L^{2l+1}_{1-\varepsilon},\Xi_{L_{1-\varepsilon}})$. 
Note that the disjoint union $(\eta_{1-\varepsilon},
\Psi_{N_{1-\varepsilon}}) \cup (f_0, \Psi_{N_0})$ is a stably
framed boundary.  

Note that the  skew-framings of
$\delta^{st}(\xi_0,\Psi_0)$ can be distingushed from the skew-framing of the immersion, which is determined by the boundary $W^{4k+3}.$ 
The two skew-framings admit
 a common Browder-Eccles invariant since they
  are well-defined as skew-framings of a common immersion
into $R^{4l+2}.$ 
Therefore for the stably framed immersion  
$(\xi_{1-\varepsilon}, \Psi_{1-\varepsilon}) \cup (\xi_0, \Psi_0),$
both invariants are trivial. 

On other hand, by the calculations,
the twisted Arf-Kervaire invariant and the twisted Browder-Eccles
invariant of $(f_0,\Psi_{N_0})$ coincide. Hence for the stably
framed immersion
 $(\eta_{1-\varepsilon},\Xi_{N_{1-\varepsilon}}),$
both invariants coincide and moreover, for the stably framed
immersion  $(f,\Xi_N) \in \Imm^{stfr}(4l+2,2l+1)^{\ast},$ both
invariants coincide. Lemma~\ref{4.1} is thus proved \qed

{\bf Proof of Theorem~\ref{main}.} Let $4l+2=30$ or $4l+2=62$.
Consider a closed framed $2l$--connected manifold
$(N^{4l+2},\Psi)$ with the Arf-Kervaire invariant 1.
Let us assume that for $l=7$  ($l=15$) there exists an embedding
$J: N^{4l+2} \subset \R^{6l+3}$  and that the normal
$(2k+2)$--bundle $\nu_{\bar I}$ of this embedding is equipped with
$9$ (resp. $10$) linearly independent sections.

The restriction of the normal bundle $\nu_J$ on an arbitrary
embedded sphere $i_{S^{2l+1}}: S^{2l+1} \to N^{4l+2}$ is the
trivial $(2l+1)$--dimensional, i.e., $15$--dimensional (resp.
$31$--dimensional) bundle (the bundle $i_{S^{2l+1}}^{\ast}(\nu_J)$
is stably dimensional and trivial), equipped with $9$ (resp. $10$)
linearly independent sections.

There is only one stably trivial but nontrivial $15$--dimensional
(resp. $31$--dimensional) bundle over $S^{15}$ (resp. $S^{31}$);
this bundle is the tangent bundle $T(S^{2l+1})$ (Kervaire and Milnor, 1963; p.
534). The bundle $i_{S^{2l+1}}^{\ast}(\nu_J)$ over $S^{15}$
($S^{31}$) is trivial if and only if 
$c(S^{15},i_{S^{15}}^{\ast}(\nu_J))=0$ (resp.
$c(S^{31},i_{S^{31}}^{\ast}(\nu_J)=0$). By the Adams theorem
(Novikov, 1964; Ch. 3, Par. 8, p.
106, the reference on the Adams result) the tangent bundle
$T(S^{15})$ (resp. $T(S^{31})$) admits no more than $8$ (resp.
$9$) linearly independent sections. By our assumption the bundle
$i_{S^{2l+1}}^{\ast}(\nu_J)$ admits $9$ (resp. $10$) linearly
independent sections. Therefore $i_{S^{2l+1}}^{\ast}(\nu_J)$ is a
trivial bundle. Hence the auxiliary homomorphism $\pi: H \to \Z/2$
is trivial and the Arf-Kervaire invariant of the stably framed
manifold $(N^{4l+2},\Xi)$ is equal to the twisted Arf-Kervaire
invariant of the pair $(J,\Xi)$, so both are equal to 1.
On the other hand, by Lemmas \ref{4.1} and \ref{4.2},  the twisted Arf-Kervaire coincides with the twisted Browder-Eccles invariant. 

Let us assume that for $l=7$  ($l=15$) there exists an embedding
$\bar J: N^{4l+2}\times I \subset \R^{6l+4}$  and that the normal
$(2k+2)$--bundle $\nu_{\bar I}$ of this embedding is equipped with
$9$ (resp. $10$) linearly independent sections over the complement
to the base segment $pt \times I \in N^{4l+2} \times I$.

The restriction of the normal bundle $\bar \nu(\bar J)$ over an
arbitrary embedded sphere $i_{S^{2l+1}}: S^{2l+1} \subset
N^{4l+2}$ is the  trivial $(2l+1)$--bundle since it is equipped with
$9$ (resp. $10$) linearly independent sections.

By the Rourke-Sanderson compression theorem (Rourke and Sanderson, 2001), we may
assume, after an appropriate isotopy of the framed embedding $\bar
J$, that the collection of segments $I$ is vertically up with
respect to the axis of projection. After the projection we obtain
an immersion $J : N^{4l+2} \looparrowright \R^{6l+3}$, which is
framed at least outside a neighborhood of a point. The double
point manifold $L^{2l+1}$ of the immersion $J$ is stably
skew-framed, hence it is in fact, framed. Let us denote this
framing by $\Xi$.

The framed manifold $(L^{2l+1},\Xi)$ determines an element of the
group $Imm^{stsf}(2l+1,2l+1)$ lying in the image of the
homomorphism $\Pi_{2l+1}=Imm^{fr}(2l+1,2l+1) \to
Imm^{stsf}(2l+1,2l+1)$. Therefore the twisted Browder-Eccles
invariant of the stably framed immersion $(J,\Psi_N)$ is equal to
the stable Hopf invariant of the framed manifold
$(L^{2l+1},\Xi_L)$. By the Toda theorem for $2l+1=15$ and by the
Adams theorem for $2l+1=31,$ the Hopf invariant is equal to $0$ (Mosher and Tangora, 1968; Sect. 18).

On the other hand, by Lemma~\ref{4.1}, the twisted Browder-Eccles
invariant of $(J,\Psi_N)$ is equal to the twisted Arf-Kervaire
invariant of $(J,\Psi_N)$. The latter is equal to the Arf-Kervaire
invariant of the framed manifold $(N,\Psi_N)$ because the
auxiliary homomorphism $\pi: H \to \Z/2$ for the stably framed
immersion $(J,\Psi_N)$ is trivial. Therefore the twisted
Arf-Kervaire invariant is equal to 1. This contradiction shows
that if the manifold $N^{4l+2} \times I$ embeds in $\R^{6l+4},$
then the collection of linearly independent sections of the normal
bundle does not exist. Theorem \ref{main} is thus proved.
 \qed
\section*{Acknowledgements}
The first author was supported by  
 the
Russian Science Foundation (project 21-11-00010). 
The second and the third author were 
supported by the Slovenian Research Agency 
grants P1-0292, J1-4031, J1-4001, N1-0278, N1-0114, and N1-0083.
We thank the referee for comments and suggestions.

\vfill\eject
\noindent
{\bf References}

\bgroup
\setlength{\parindent}{-1em} 

\vskip 0.5cm

\begin{description}
\item
Adams, J. F. (1962), "Vector fields on spheres", {\it Annals of Mathematics} (2)
{\bf 75}:3 , 603--632.
URL: \url{http://www.jstor.org/stable/1970213?origin=JSTOR-pdf}

\item
Akhmetiev, P.~M. (2016),
"Kervaire problem in stable homotopy theory", 
  arXiv:1608.01206v1 [math.AT].
  DOI \url{ https://doi.org/10.48550/arXiv.1608.01206}

\item
Akhmetiev, P.~M. (2022),
"Geometric approach to stable homotopy groups of spheres II; Arf-Kervaire Invariants", 
arXiv:2110.04542v3 [math.AT], DOI: \url{https://doi.org/10.48550/arXiv.2110.04542}

\item
Akhmetiev, P.~M., Cencelj, M., and Repov\v{s}, D.~D. (2010), 
"The Arf-Kervaire invariant of framed manifolds as an obstruction to embeddability", 
	arXiv:0804.3164 [math.AT].
	DOI: \url{https://doi.org/10.48550/arXiv.0804.3164}
	
\item
	Akhmetiev, P.~M. and Eccles, P.~J., (1999), "A geometrical proof of
Browder's result on the vanishing of the Kervaire invariant",
{\it Proceedings of the Steklov Institute of Mathematics} {225}, 40--44.
URL: \url{https://www.researchgate.net/publication/2619825_A_geometrical_proof_of_Browder
\%27s_result_on_the_vanishing_of_the_Kervaire_invariant}

\item
Akhmetiev, P.~M. and Eccles, P.~J., (2007) , "The relationship between
framed bordism and skew-framed bordism", {\it Bulletin of the London Mathematical Society}
{\bf 39} , 473--481.
URL: \url{https://doi.org/10.1112/blms/bdm031}

\item
Arnold, V.~I., Varchenko, A.N., and Gusein-Zade, S.M. (1982), 
{\it  Singularities of Differentiable Maps: Classification of Critical Points, Caustics and Wave Fronts}, Moscow:
Nauka. (in Russian)
URL: \url{https://doi.org/10.1007/978-0-8176-8340-5}

\item
Browder, W. (1969), "The Kervaire invariant of framed manifolds and
its generalization", {\it Annals of Mathematics} (2) 90, 157--186.
DOI: \url{https://doi.org/10.2307/1970686}

\item
Browder, W. (1972), {\it Surgery on Simply Connected
Manifolds}, Ergebnisse der Mathematik und ihrer Grenzgebiete, vol. 65. Berlin:Springer
DOI: \url{https://doi.org/10.1007/978-3-642-50020-6}

\item
 Eccles, P.~J. (1979),  "Filtering framed bordism by embedding
codimension", {\it Journal of the London Mathematical Society} (2) 19, 163--169.
DOI: \url{https://doi.org/10.1112/jlms/s2-19.1.163}
	
\item
Eccles, P.~J.  (1981),  "Codimension one immersions and the Kervaire
invariant one problem", {\it Mathematical Proceedings of the Cambridge Philosophical Society} 90, 483--493.
DOI: \url{https://doi.org/10.1017/S0305004100058941}	

\item
Haefliger, A. (1962), "Plongements diff\'erentiables de vari\'et\'es
dans vari\'et\'es", {\it Commentarii Mathematici Helvetici} 36, 47--82.
DOI: \url{https://doi.org/10.1007/BF02566892}

\item
Hill, M.~A., Hopkins, M.~J., and Ravenel, D.~C. (2016),
 "On the non-existence of elements of Kervaire invariant one", 
{\it Annals of Mathematics} (2) 184, 1-262.
DOI: \url{https://doi.org/10.4007/annals.2016.184.1.1}

\item
Hirsch, M.~W. (1961), "The imbedding of bounding manifolds in Euclidean
space",  {\it  Annals of Mathematics} (2) 74, 494--497.
DOI: \url{https://doi.org/10.2307/1970293}

\item
Jones, J. and Rees, E. (1978), "Kervaire's invariant for framed
manifolds", In: {\it Algebraic and Geometric Topology, Proceedings of
Symposia in Pure Mathematics}, Vol. 32, Ed. by R. James Milgram, Providence RI: American Mathematical Society,
  141--147.
URL: \url{https://mathscinet.ams.org/mathscinet-getitem?mr=520501}

\item
Kervaire, M. and Milnor, J. (1963), "Groups of homotopy spheres I",
{\it Annals of Mathematics} (2) 77, 504--537.
DOI: \url{https://doi.org/10.2307/1970128}

\item
 Kreck, M. (2000), "A guide to the classification of manifolds", In:
{\it Surveys on Surgery Theory, Vol. 1}, Annals of Mathematical Studies, Vol. 145, Ed. by S. Cappell, A. Ranicki, and J. Rosenberg,
Princeton, NJ: Princeton Univiversity Press,  121--134.
URL: \url{https://mathscinet.ams.org/mathscinet-getitem?mr=1747533}

\item
Lin, W.-H.  (2001),  "A proof of the strong Kervaire invariant in
dimension 62", In: {\it First International Congress of Chinese
Mathematicians (Beijing, 1998), AMS/IP Studies in Advanced Mathematics Vol. 20},
Providence, RI: American Mathematical Society,   351--358.
URL: \url{https://mathscinet.ams.org/mathscinet-getitem?mr=1830191}

\item
 Mosher, R.~S. and Tangora, M.~C.  (1968), {\it Cohomology Operations and
Applications in Homotopy Theory}, New York: Harper-Row Publishers.
URL: \url{https://mathscinet.ams.org/mathscinet-getitem?mr=226634}

\item
Novikov, S. P.  (1964), "Homotopy equivalent smooth manifolds I", {\it Izvestiia Akademii Nauk SSSR, Seriya Matematicheskaya} 28, 365--474.
 (in Russian) 
English translation: URL: \url{https://homepage.mi-ras.ru/~snovikov/10.pdf}

\item
Pontryagin, L.~S. (1955), {\it Smooth manifolds and their applications in
homotopy theory}, Providence, RI:  American Mathematical Society Translations (Ser. 2)  11
1--114.
\url{https://mathscinet.ams.org/mathscinet-getitem?mr=115178},
\url{https://mathscinet.ams.org/mathscinet-getitem?mr=71767}

\item
Randall, D. (1999), "Embedding homotopy spheres and the Kervaire
invariant", In: {\it Homotopy Invariant Algebraic Structures (Baltimore,
MD, 1998), Contemporary Mathematics} Vol. 239, Providence,
RI: American Mathematical Society, 239--243.
URL: \url{https://mathscinet.ams.org/mathscinet-getitem?mr=1718084}

\item
Rourke, C. P.  and Sanderson, B. J.  (2001), "The compression theorem,
I.", {\it Geometry and Topology} 5, 399--429.
DOI: \url{http://dx.doi.org/10.2140/gt.2001.5.399}	
	
\item
Wells, R. (1966), "Cobordism groups of immersions", {\it  Topology} 5,
281--294.
DOI: \url{https://doi.org/10.1016/0040-9383(66)90011-5}
	
\end{description}
 
\vskip 1.0cm

\noindent
Tikhonov Moscow Institute of Electronics and Mathematics (MIEM HSE),
Moscow, 109028
\& 
Pushkov Institute of Terrestrial Magnetism, Ionosphere and Radio Wave Propagation of the
Russian Academy of Science (IZMIRAN),
 Moscow region, Troitsk, 108840,
Russian Federation
\\
{\sl Email address}: {pmakhmet@mail.ru }

\vskip 0.5cm

\noindent  
Faculty of Education and Faculty of Mathematics and Physics, University of Ljubljana,  1000,
\&
Institute of Mathematics, Physics and Mechanics,
Ljubljana 1001, Slovenia\\
{\sl Email address}: {matija.cencelj@guest.arnes.si}

\vskip 0.5cm

\noindent 
Faculty of Education and Faculty of Mathematics and Physics, University of Ljubljana,  1000,
\&
Institute of Mathematics, Physics and Mechanics,
Ljubljana 1001, Slovenia\\
{\sl Email address}: {dusan.repovs@guest.arnes.si}
\end{document}